\documentclass[11pt]{amsproc}
\usepackage{amssymb,longtable,amscd,euscript}
\newcommand{\xref}[1]{{\rm \ref{#1}}}

\newcommand{\muu}{\mbox{\boldmath $\mu$}}
\newcommand{\stackunder}[2]{\mathrel{\mathop{#2}\limits_{#1}}}
\newcommand{\mt}[1]{\operatorname{#1}}

\newcommand{\OOO}{\mathcal{O}}
\newcommand{\PP}{{\mathbb P}}
\newcommand{\CC}{{\mathbb C}}
\newcommand{\QQ}{{\mathbb Q}}
\newcommand{\NN}{{\mathbb N}}
\newcommand{\KKK}{{\EuScript{K}}}

\newcommand{\ZZ}{{\mathbb Z}}
\newcommand{\var}{{\varphi}}
\newcommand{\ord}{{\operatorname{ord}}}
\newcommand{\mult}{{\operatorname{mult}}}
\newcommand{\Diff}{{\operatorname{Diff}}}
\newcommand{\Center}{{\operatorname{Center}}}
\newcommand{\discr}[1]
{{\operatorname{discr}}\left(#1\right)}
\newcommand{\lin}{\text{---}}
\newcommand{\wt}{\mathrm{wt}}

\newcommand{\y}[1]{#1}

\newcommand{\comment}[1]{}
\newtheorem{theorem}[subsection]{Theorem}
\newtheorem{proposition}[subsection]{Proposition}
\newtheorem{lemma}[subsection]{Lemma}
\newtheorem{corollary}[subsection]{Corollary}

\theoremstyle{definition}
\newtheorem*{definition*}{Definition}
\newtheorem{definition}[subsection]{Definition}
\newtheorem{example}[subsection]{Example}
\theoremstyle{remark} 

\newtheorem{remarks}[subsection]{Remarks}
\renewcommand\labelenumi{(\roman{enumi})}

\author{Yuri~G.~Prokhorov}
\title{Elliptic Gorenstein singularities, log canonical thresholds 
and log Enriques surfaces}
\thanks{The author was
partially supported by grants INTAS-OPEN-97-2072
and RFFI 99-01-01132}

\address{
Department of Algebra, Faculty of Mathematics, Moscow State
Lomonosov University, Moscow 117234, Russia}
\email{prokhoro@mech.math.msu.su}
\address{Current Address: Max-Plank-Institut f\"ur Mathematik,
Vivatsgasse 7, 53111 Bonn, Germany}
\email{prokhoro@mpim-bonn.mpg.de}

\date{}
\begin{document}
\begin{abstract}
We point out an 
interesting relation between 
hypersurface elliptic singularities and log Enriques surfaces:
with a few exceptions, every hypersurface elliptic singularity
define some klt log Enriques surface $(S,\Diff)$. In many cases, 
the log canonical cover of $(S,\Diff)$ is a Du Val K3 surface with
``large'' nonsymplectic automorphism group.
\end{abstract}
\maketitle 
\section{Introduction}
The aim of this paper is to point out an 
interesting relation between 
hypersurface elliptic singularities and log Enriques surfaces:
given hypersurface elliptic singularity $F\subset \CC^3$,
with a few exceptions, 
there is a unique exceptional divisor $S$ over $\CC^3$
with discrepancy $a(S,cF)=-1$, where $c:=c(\CC^3,F)$ is the log canonical 
threshold. Further, there is a uniquely defined blowup
$\var\colon Y\to \CC^3$ of $S$ such that the pair $(S,\Diff_S(cF_Y))$
is a klt log Enriques surface \cite{Z3}, \cite{B}.

Our main results are the following two theorems.
For basic definitions we refer to Section~\xref{sect-p}.
\begin{theorem}
\label{main1}
Let $0\in F\subset \CC^3$ a hypersurface elliptic singularity
given by the equation $f(x,y,z)=0$ and 
let $c:=c(F)$ be the log canonical threshold. 
Then one of the following holds:
\begin{enumerate}
\item
$c=1$ and the singularity $(F\ni 0)$ is log canonical;
\item
$c=5/6$, $\mult_0 f(x,y,z)=3$, $f_3(x,y,z)$ has a triple factor
and the minimal resolution has one of the following five 
dual graphs (see \cite[Table~3]{L}):
\[
2\mathrm{A}_{1*o}\mathrm{A}_{4*o},\quad
2\mathrm{A}_{1*o}\mathrm{E}_{6o},\quad
2\mathrm{A}_{4*o},\quad
\mathrm{A}_{4*o}\mathrm{E}_{6o},\quad
2\mathrm{E}_{6o};
\]
\item
the pair $(\CC^3,cF)$ is exceptional, that is there is exactly one 
exceptional divisor $S$ with discrepancy $a(S,cF)=-1$.
\end{enumerate}
\end{theorem}

The most interesting is case (iii):

\begin{theorem}
\label{main2}
Let $0\in F\subset \CC^3$ a hypersurface elliptic singularity
and let $c:=c(F)$ be the log canonical threshold. 
Assume that $(F\ni 0)$ is such as in {\rm (iii)} of
Theorem~\xref{main1}.
\begin{enumerate}
\item
There is a (unique) blowup $\var\colon Y\to \CC^3$ with 
$\rho(Y/\CC^3)=1$, $\QQ$-factorial $Y$ and exceptional divisor 
$S$, where $a(S,cF)=-1$.
\item
Let $F_Y$ be the proper transform of $F$ on $Y$. Then
$(Y,S+cF_Y)$ is plt.
\item
$(S,\Diff_S(cF_Y))$ is a klt log Enriques surface that is
$S$ is a normal projective surface, $\Diff_S(cF_Y)$ is a boundary on $S$
such that $(S,\Diff_S(cF_Y))$ is klt and $K_S+\Diff_S(cF_Y))$ is 
numerically trivial.
\item
Furthermore, assume that $c$ is of \emph{standard form}, i.e., $c=1-1/m$,
for $m\in \NN$. Then there is the following commutative diagram
\begin{equation}
\label{cd}
\begin{CD}
Y'@>\psi>>Y\\
@V{\var'}VV@V{\var}VV\\
X'@>\pi>>\CC^3\\
\end{CD}
\end{equation}
where $X'$ is a hypersurface simple $K3$-singularity (see 
Definition~\xref{def-K3}),
$\pi$ is a cyclic $\muu_m$-covering ramified along $F$,
$\var'$ is the blowup of a unique divisor $S'$ with discrepancy 
$a(S',0)=-1$. The pair $(Y',S')$ is canonical of index $1$ and
$Y'$ has terminal singularities.
Moreover, $\psi(S')=S$, $S'$ is a normal $K3$-surface with 
Du Val singularities and $S'\to S$ is a ramified
$\muu_m$-covering.
\end{enumerate}
\end{theorem}

\begin{remarks}
\begin{enumerate}
\item
Surfaces $S'$ gives a lot of 
examples of $K3$-singularities with ``large''
non-symplectic automorphism groups. 
\item
In many cases we have $c>6/7$. 
Then $(S,\Diff_S(c-\varepsilon)F_Y)$ for $0< \varepsilon\ll1$
is an \emph{exceptional} log del Pezzo surface \cite{Sh1},
\cite{KeM}.\footnote{The notion of log del Pezzo surface
was introduced in \cite{MT}.} 
\item
Note that replacing in Theorem~\xref{main1}
``elliptic'' with ``rational'' we get the non-interesting case
$c=1$.
\item
Three-dimensional hypersurface simple K3 singularities 
given by a non-degenerate function were classified by
Yonemura \cite{Y}. The list contains exactly 95 families 
which are related to the famous 95 families of
weighted K3 hypersurfaces of Fletcher and Reid \cite{Flet},
\cite[Sect.~4.5]{Reid-canonical}. 
In Sect. \xref{tabl} we give the correspondence between 
hypersurface elliptic singularities with standard $c$ and 
the list in \cite{Y}. Note, however, that not all
95 families can be obtained from our construction.
\item
A construction similar to that of
Theorem~\xref{main2} is possible for all singularities 
$F\subset\CC^3$ under the assumption that the log canonical threshold
$c(\CC^3,F)$ does not come from a two-dimensional log 
singularity (cf. \cite[Prop. 1.3]{Pr-n}). 
\end{enumerate}
\end{remarks}

The paper is organized as follows. Section \xref{sect-p}
is auxiliary. Section~\xref{proof-T1} is most important.
It contains estimates of log canonical thresholds which gives
the proof of Theorem~\xref{main1}. Theorem~\xref{main2} is 
an easy consequence of Theorem~\xref{main1}. It is proved in 
Section~\xref{proof-T2}. Finally, in Section~\xref{tabl},
we present tables containing concrete computations
for our construction
for all sample equations of elliptic hypersurface singularities 
listed in \cite{L}.

\subsubsection*{Acknowledgments}
This work was carried out during my stay at 
Max-Planck-Institut f\"ur Mathematik in Bonn.
I would like to thank MPIM for wonderful 
working environment.
I also would like to thank Professor De-Qi Zhang
for pointing me out some references.

\section{Preliminaries}
\label{sect-p}
All varieties are assumed to be algebraic and defined over $\CC$.
Notation of the Log Minimal Model Program
are used freely.
\begin{definition}[\cite{R0}, \cite{L}]
Let $(F\ni o)$ be a normal surface singularity and
let $\mu\colon \tilde F\to F$ be the minimal resolution.
$(F\ni o)$ is said to be \emph{elliptic} if 
$R^1\mu_*\OOO_{\tilde F}\simeq \CC$.
\end{definition}

\begin{definition}
\label{def-K3}
A normal log canonical singularity $(X\ni o)$ is said to be 
\emph{exceptional} if for any boundary $D$ such that $(X,D)$ is
lc there is at most one divisor $S$ of the function field 
$\KKK(X)$ with discrepancy $a(E,D)=-1$. An isolated log canonical
non-klt 
Gorenstein 
exceptional singularity is called \emph{simple $K3$ singularity}
\cite{IW}\footnote{The definition given in \cite{IW}
uses the plurigenera and types of Hodge structures
of essential divisors. By the main result of \cite{IW}
this is equivalent to our definition}.
\end{definition}

\begin{definition}
Let $X$ be a normal variety with at
worst log canonical singularities and let $F$ be an effective
non-zero 
$\QQ$-Cartier divisor on $X$. \textit{The log canonical
threshold} of $(X,F)$ is defined by
\[
c(X,F)=\sup \left\{c \mid \text{$(X,cF)$
is log canonical}\right\}.
\]
It is easy to see that $c(X,F)$ is a rational number 
(see \cite{K}). Moreover, $c(X,F)\in[0,1]$ 
whenever $F$ is integral. 
We frequently write $c(F)$ instead of $c(X,F)$ if no confusion is likely.
\end{definition}

\begin{example}
\label{ex}
Let $C\subset \CC^2$ be a plane curve given by
$x^n=y^2$. Then
\[
c(\CC^2,C)=\frac12+\frac1n.
\]
Indeed, consider the case when $n$ is even: $n=2m$.
Then $C=C_1+C_2$, where $C_1$ and $C_2$ are non-singular curves 
given by $y\pm x^m=0$. The minimal good resolution of 
$C\subset \CC^2$ is the sequence of $m$ blowups having the 
following dual graph:
\[
\begin{array}{ccccccc}
\stackrel{E_1}\circ&\lin&\stackrel{E_2}\circ&
\lin\cdots\lin&\stackrel{E_m}\circ&\lin&\stackrel{C_1}\bullet\\
&&&&|&&\\
&&&&\stackunder{C_2}\bullet&&\\
\end{array}
\]
It is easy to show by induction 
that $a(E_i,tC)=-2ti+i=i(1-2t)$.
Hence
\[
c(\CC^2,C)=\sup_t\left\{ t\mid \min_i\{a(E_i,tC)\}\ge -1\right\}=
\sup_t\Bigl\{t \mid m(1-2t)\ge -1\Bigr\}=\frac12+\frac1{2m}.
\]
The computations in the case when $n$ is odd are similar.
\end{example}

\begin{proposition}[\cite{R0}, \cite{L}, see also \cite{R}, \cite{KM}]
Let $(F\ni o)$ be a Gorenstein elliptic surface singularity and
let $\mu\colon \tilde F\to F$ be the minimal resolution.
Write $K_{\tilde F}=\mu^*K_F-\sum z_iA_i$, where the $A_i$ are prime 
exceptional divisors and the $z_i$ are non-negative integers.
Then $\sum z_iA_i$ coincides with the fundamental cycle $Z$ of
$(F\ni o)$.
\end{proposition}

\begin{theorem}[\cite{R0}, \cite{L} see also \cite{KM}]
\label{cla}
Let $(F\ni o)$ be a Gorenstein elliptic surface singularity and
let $\mu\colon \tilde F\to F$ be the minimal resolution. Let 
let $Z=\sum z_iA_i$ be the fundamental cycle and let $d:=-Z^2$. 
\begin{enumerate}
\renewcommand\labelenumi{(\roman{enumi})}
\item
Assume $d\ge 3$. Then $(F\ni o)$ has multiplicity $d$
and embedding dimension $d$. For any embedding 
$(F\ni o)\hookrightarrow (\CC^{d}\ni 0)$ consider the weight 
$w$ with $w(x_1,\dots,x_n)=(1,\dots,1)$.
\item
Assume $d= 2$. Then $(F\ni o)$ has multiplicity $2$
and embedding dimension $3$. After an analytic coordinate 
change it can be given by an equation
\[
z^2+q(x,y)=0,\quad\text{where}\quad \mult_0q\ge 4.
\] 
Let $w$ be the weight $w(x,y,z)=(1,1,2)$.
\item
Assume $d= 1$. Then $(F\ni o)$ has multiplicity $2$
and embedding dimension $3$. After an analytic coordinate 
change it can be given by an equation
\[
z^2+y^3+yq(x)+r(x)=0,\quad\text{where}\quad \mult_0q\ge 4,\ 
\mult_0r\ge 6.
\] 
Let $w$ be the weight $w(x,y,z)=(1,2,3)$.
\end{enumerate}
Let $\sigma\colon \bar F\to F$ be the weighted blowup 
(see e.g. \cite[4.56]{KM})
with weight 
$w$. Then $\bar F$ has only Du Val singularities,
it is dominated by $\tilde F$ via $\tilde\mu\colon \tilde F\to \bar F$
and $K_{\bar F}=-\tilde\mu_*Z$. 
\end{theorem}

\section{Proof of Theorem~\ref{main1}}
\label{proof-T1}
We fix the following notation. 
\subsection{Notation}
\label{not}
Let $\left(0\in F\subset \CC^3\right)$ be 
a hypersurface elliptic singularity
with invariant $d:=-Z^2$, $1\le d\le 3$ 
and let $c:=c(F)$ be the log canonical threshold. 
Let $f(x,y,z)=0$ be an equation defining $F$.
We denote the degree $m$ homogeneous part of a polynomial 
$g(x,y,z)$ by $g_m(x,y,z)$.

\begin{lemma}
The following are equivalent:
\begin{enumerate}
\item
the singularity $(F\ni o)$ is log canonical;
\item
$c=1$.
\end{enumerate}
\end{lemma}
\begin{proof}
Follows by \cite[Cor. 17.12]{Ut}
\end{proof}

From now on we assume that
$c<1$ and $(F\ni o)$ is not log canonical.
The key point in the proof 
of Theorem~\xref{main1} is to show that 
in many cases $c\neq 5/6$. Then we can use the following 
weaker version of \cite[Prop.~1.3]{Pr-n}.

\begin{proposition}[{\cite{Pr-n}}]
\label{main12}
Let $(X\ni o)$ be a three-dimensional $\QQ$-factorial
log terminal singularity and let 
$F$ be an (integral) Weil divisor on $X$.
Assume that $c:=c(X,F)>1/2$.
Then one of the following holds:
\begin{enumerate}
\item
$c=\frac12+\frac1n$, $n\in\ZZ$, $n\ge 2$; or
\item
$c\neq\frac12+\frac1n$, $n\in\ZZ$ and there is exactly one divisor 
$S$ of the function field 
$\KKK(X)$ with discrepancy $a(S,cF)=-1$ (i.e., the pair 
$(X,cF)$ is \emph{exceptional} in the sense of \cite{Sh}).
\end{enumerate}
Moreover in case {\rm (ii)}, $\Center(S)=o$. 
\end{proposition}

The following proposition gives a rough estimate
for $c(F)$ in our case. 

\begin{proposition}
\label{prop-est}
Notation as in \xref{not}. 
\begin{enumerate}
\item
If $d=-1$, then $c>11/12$.
\item
If $d=-2$, then $c>5/6$.
\item
Assume that $d=-3$. 
Then $c>7/9$. Moreover,
\begin{enumerate}
\item
If $f_3(x,y,z)$ has no multiple factors, then $c>5/6$.
\item
If $f_3(x,y,z)$ has no triple factors, then $c>4/5$.
\end{enumerate}
\end{enumerate}
\end{proposition}

\begin{proof}
Let
\[
\begin{array}{rccc}
\sigma\colon W&\longrightarrow&\CC^3\\
\cup&&\cup\\
\sigma_F\colon\bar F&\longrightarrow&F\\
\end{array}
\]
be the blowup as in Theorem \xref{cla}, where 
$\bar F$ is the proper transform of $F$. 
Let $E$ be the exceptional divisor.
Note that $W$ has at worst isolated singularities
(because the weights $w_i$ in Theorem \xref{cla} are prime to each other). 
Hence $\Diff_E(0)=0$. 
Clearly, $E$ is the weighted projective plane 
(see e.g. \cite{Flet})
$\PP(1,2,3)$, $\PP(1,1,2)$ and $\PP^2$
in cases $d=1$, $2$ and $3$, respectively.
Put $C:=\bar F|_E$. 
In all cases, $a(E,F)=-1$. Therefore,
$\sigma^*(K_{\CC^3}+F)=K_W+\bar F+E$, so 
$K_E+C\sim 0$. In particular, $p_a(C)=1$.
By Theorem \xref{cla} and the Inversion of Adjunction \cite[17.6]{Ut}
the pair $(W,\bar F)$ is plt.

(i)
In this case $C$ is given by the equation $z^2+y^3+ayx^4+bx^6=0$, $a,b\in\CC$.
It is easy to see that $C$ is reduced, irreducible and contained into 
the smooth locus of $E$. 
Since $p_a(C)=1$, the curve $C$ have at
most one singular point which is an ordinary double point or a 
simple cusp. In both cases $(E,\frac56C)$ is lc. By the Inversion of 
Adjunction \cite[17.7]{Ut}, so is $(W,E+\frac56\bar F)$.
Using convexity of the lc property \cite[1.3.2]{Sh} one can show that
the pair
\begin{equation}
\label{eq-pir}
\left(W,\frac12\left(E+\frac56\bar F\right)+\frac12\bar F\right)=
\left(W,\frac12E+\frac{11}{12}\bar F\right)
\end{equation}
is klt.
Standard toric technique gives us
\[
a\left(E,\frac{11}{12}F\right)=
-1+1+2+3-\frac{11}{12}\cdot w\mt{-ord}(f)=-\frac12.
\]
Hence we can write
\[
K_W+\frac12 E+\frac{11}{12}\bar F=
\sigma^*\left(K_{\CC^3}+\frac{11}{12}F\right).
\] 
Thus 
$(\CC^3,\frac{11}{12} F)$ is klt \cite[3.10]{K}.

(ii)
Here $C$ is given by the equation $z^2+q_4(x,y)=0$.
We claim that $(E,\frac12C)$ is lc. 
As above, $C$ is contained into 
the smooth locus of $E$ and $p_a(C)=1$.
If $C$ is reduced and irreducible, we have that 
$(E,\frac56C)$ is lc (because in this case
the singularities of $C$ are at worst ordinary double points or
a simple cusp, cf. Example~\xref{ex}). 
Assume that $C$ is reduced but not irreducible.
Then $C=C_1+C_2$, where $C_1$ and $C_2$ are smooth rational curves
and $C_1\cdot C_2=2$. Hence $c(E,C_1+C_2)\ge 3/4$ (see Example~\xref{ex}). 
Finally,
if $C$ is reduced, then again $C=2C_1$ and $C_1$ is smooth.
Hence $(E,\frac12C)=(E,C_1)$ is plt. This proves our claim.
By Inversion of Adjunction \cite[17.7]{Ut}, $(W,E+\frac12\bar F)$ is lc. 
Similar to \eqref{eq-pir} the pair
\[
\left(W,\frac13\left(E+\frac12\bar F\right)+\frac23\bar F\right)=
\left(W,\frac13E+\frac56\bar F\right)
\] 
is klt. So, 
\[
K_W+\frac13E+\frac56\bar F=\sigma^*\left(K_{\CC^3}+\frac56F\right)
\]
is also klt.

(iii)
Here $C$ is a cubic curve given by the equation $f_3(x,y,z)=0$.
Let $c':=c(C,E)$. 
We have $c'\ge 1/3$ (see e.g. \cite[Lemma~8.10]{K}).
Moreover, $c'\ge 2/3$ whenever $C$ is reduced (cf. Example~\xref{ex})
and $c'\ge 1/2$ whenever $C$ has no components of multiplicity 
$3$. Similar to \eqref{eq-pir} the pair
\[
\left(W,
\frac{1}{4-3c'}\left(E+c'\bar F\right)+
\frac{3-3c'}{4-3c'}\bar F\right)=
\left(W,\frac{1}{4-3c'}E+
\frac{3-2c'}{4-3c'}\bar F\right)
\] 
is klt. Put $t:=(3-2c')/(4-3c')$. Then $a(E,tF)=2-3t$.
Hence 
\[
K_W+(3t-2)E+t\bar F=\sigma^*(K_{\CC^3}+tF).
\] 
Since 
\[
3t-2\le \frac{1}{4-3c'}
\qquad \text{and}\qquad
t\le \frac{3-2c'}{4-3c'},
\]
both pairs 
$(W,(3t-2)E+t\bar F)$ and $(\CC^3,tF)$ are klt.
This proves the proposition.
\end{proof}

The following is a consequence of Propositions~\xref{main12}
and \xref{prop-est}.

\begin{corollary}
\label{cor}
Notation as in \xref{not}. Assume that $(X,cF)$ is 
not exceptional. 
Then $Z^2=-3$, $c=5/6$ and $f_3(x,y,z)$ has a multiple 
factor.
\end{corollary}

Thus in cases $d=1$ and $d=2$ Theorem~\xref{main1} 
follows by Corollary~\xref{cor}. It remains to consider the case
$d=3$ with $f_3$ having multiple factors.
This will be done in \xref{3ple} and \xref{2ple}.

\subsection{Cases: $f_3$ has a triple factor}
\label{3ple}
We may assume that $f_3(x,y,z)=z^3$.
Consider the $w$-blowup 
$\varphi\colon Y\to \CC^3$ with $w:=(3,3,4)$.
Let $S$ be the exceptional divisor and
let $F_Y$ be the proper transform of $F$.
It is easy to see that $\ord_wf(x,y,z)=12$ and
\[
f_{w-\wt=12}(x,y,z)=z^3+f_4(x,y),
\]
where is the weighted multiplicity of $f$ at $0$
and $f_{w-\wt=12}$ is the weighted homogeneous 
term of weighted degree $12$.
Therefore,
\[
a\left(S,\frac56F\right)=-1+3+3+4-\frac56\cdot12=-1.
\]
In particular, $c\le 5/6$. Assume that $(\CC^3,cF)$ is not exceptional.
Then by Corollary~\xref{cor}, we have $c=5/6$ and 
$\left(S,\Diff_{S}(\frac56F_Y)\right)$ is lc but not klt.
Further, $S\simeq\PP(3,3,4)$. 
Up to isomorphism we may assume that
the weighted projective plane $S$ is written in the 
\emph{well-formed} form 
(see e.g. \cite{Flet}):
$S\simeq\PP(1,1,4)$. In this case,
$\Diff_{S}(\frac56F_Y)=\frac56L+\frac23L_3$, where $L_3:=\{z=0\}$
and $L:=\{z+f_4(x,y)=0\}$. 
It is easy to check (cf. Example~\xref{ex}) 
that $(S,\frac56L+\frac23L_3)$ is 
klt (resp. lc) if and only if $L\cup L_3$ have 
only normal crossings (resp. $L_3$ and $L$ meet each other either
normally or tangentially to the first order). Therefore,
$f_4(x,y)$ has a double factor but has no triple factors.
Thus we may assume that $f_4(x,y)=y^4+x^2y^2$ or $f_4(x,y)=x^2y^2$. So,
\[
f(x,y,z)=
\begin{cases}
\text{a)}\quad z^3+y^4+x^2y^2+g(x,y,z),&\text{or}\\
\text{b)}\quad z^3+x^2y^2+g(x,y,z),&\\
\end{cases}
\] 
where $g(x,y,z)=(\text{terms of $w$-degree $>12$})$.
Now consider the blowup $\sigma\colon W\to \CC^3$. In the affine chart 
over $(x\neq 0)\subset \CC^3$ the surface 
$\bar F\subset W_{x\neq 0}\simeq \CC^3$ is given by
\[
\bar f(x,y,z)=
\begin{cases}
\text{a)}\quad z^3+xy^4+xy^2+x^{-3}g(x,xy,xz),&\text{or}\\
\text{b)}\quad z^3+xy^2+x^{-3}g(x,xy,xz).&\\
\end{cases}
\] 
On the other hand, $\bar F$ has only Du Val singularities.
Hence $\bar f$ must contain a term of degree $2$. 
It can be either $xz$ or $x^2$. So,
\[
g(x,y,z)=
\begin{cases}
x^3z+\cdots,&\text{or}\\
x^5+\cdots.&\\
\end{cases}
\] 
In case b), considering the chart $(y\neq 0)$ we get that
$g(x,y,z)$ also must contain either $y^3z$ or $y^5$.
Finally,
\begin{equation}
\label{eq-mnogo}
f(x,y,z)=
\begin{cases}
z^3+y^4+x^2y^2+x^3z+\cdots,&\\
z^3+y^4+x^2y^2+x^5+\cdots,&\\
z^3+x^2y^2+x^3z+y^3z+\cdots,&\\
z^3+x^2y^2+x^5+y^3z+\cdots,&\text{or}\\
z^3+x^2y^2+x^5+y^5+\cdots.&\\
\end{cases}
\end{equation} 
These are cases $2\mathrm{A}_{1*o}\mathrm{A}_{4*o}$,
$2\mathrm{A}_{1*o}\mathrm{E}_{6o}$,
$2\mathrm{A}_{4*o}$,
$\mathrm{A}_{4*o}\mathrm{E}_{6o}$, and
$2\mathrm{E}_{6o}$, respectively (see the last section
of Table~3 in \cite{L}). Conversely, if $f(x,y,z)$
is one of the equations \eqref{eq-mnogo}, then $f_4(x,y)$
has a double, not triple factor. By the above,
$\left(S,\Diff_{S}(\frac56F_Y)\right)$ is lc but not klt.
So, the pair $(X,cF)$ is not exceptional in these cases.

\subsection{Cases: $f_3$ has a double factor
and a single factor}
\label{2ple}
Then $E\cap \bar F$ is given (scheme-theoretically) 
by $f_3(x,y,z)=0$ on
$E\simeq\PP^2$. Hence $E\cap \bar F=2C_1+C_2$,
where $C_1$ and $C_2$ are lines on $\PP^2$. 
Since $a(E,\frac56F)=2-\frac56\cdot \mt{ord}_0 f=\frac12$,
we can write
\begin{equation}
K_{W}+\frac56\bar F+\frac12E=\sigma^*\left(K_X+\frac56F\right).
\end{equation}
Therefore, $c\ge 5/6$ (resp. $c>5/6$) if and only if 
$K_{W}+\frac56\bar F+\frac12E$ is lc (resp. klt).

\begin{lemma}
\label{l-exc}
If $(\bar F,\frac12C_1+\frac14C_2)$ is lc, then 
the pair $(X,cF)$ is exceptional.
\end{lemma}

\begin{proof}
By the Inversion of Adjunction \cite[17.6, 17.7]{Ut} we have
that $\left(W,\frac14E+\bar F\right)$ is lc.
Assume that the pair $(X,cF)$ is nonexceptional.
By Corollary~\xref{cor}, $c=5/6$.
Hence, there are at least two divisors $E_1$ and $E_2$ such that 
$a\left(E_i,\frac12E+\frac56\bar F\right)=
a\left(E_i,\frac56 F\right)=-1$.
Note that
\begin{equation}
\label{combination}
\frac12E+\frac56\bar F=
\frac13\left(E+\frac12\bar F\right)+
\frac23\left(\frac14E+\bar F\right).
\end{equation}
Since both pairs $\left(W,E+\frac12\bar F\right)$ and 
$\left(W,\frac14E+\bar F\right)$ are lc, we have 
\[
a\left(E_i,E+\frac12\bar F\right)\ge -1
\quad\text{and}\quad
a\left(E_i,\frac14E+\bar F\right)\ge -1.
\]
Taking into account that the discrepancy
\[
a\left(E_i,t\left(E+\frac12\bar F\right)+
(1-t)\left(\frac14E+\bar F\right)\right)
\]
is a linear function in $t$,
we obtain 
\begin{equation*}
\label{eq-a}
a\left(E_i,\frac12E+\frac56\bar F\right)=
a\left(E_i,E+\frac12\bar F\right)=
a\left(E_i,\frac14E+\bar F\right)=-1.
\end{equation*}
Clearly, $\Center_W(E_i)\subset C_1$.
Again by the Inversion of Adjunction, 
$\left(W,\frac14E+\bar F\right)$ is plt in codimension 
two. If $\Center_W(E_i)=C_1$, then 
$a\left(E_i,\frac14E+\bar F\right)>-1$.
Hence $a\left(E_i,\frac12E+\frac56\bar F\right)>-1$, a contradiction.
Therefore, $\Center_W(E_i)$ is a point, say $P_i$.
So, $\discr{\Center\subset P_i,W,E+\frac12\bar F}=-1$.
By \cite[Th. 17.2, Cor. 17.11]{Ut}, 
\[
\discr{\Center\subset P_i,E,\Diff_E\left(\frac12\bar F\right)}=-1.
\]
On the other hand, $\left(E,\Diff_E(\frac12\bar F)\right)=
\left(\PP^2,C_1+\frac12C_2\right)$ is plt. 
The contradiction proves our lemma.
\end{proof}

\subsubsection*{}
Now we prove that in many cases 
$(\bar F,\frac12C_1+\frac14C_2)$ is lc:

\begin{lemma}
\label{C1C2}
Notation as above.
The pair 
$(\bar F,\frac12C_1+\frac14C_2)$ is lc
in the following cases (see Sect.~{\rm 4} of
Table~{\rm 2} in \cite[pp. 1293--1294]{L}):
\begin{multline}
\label{mul-cacses}
\text{
\begin{tabular}{p{300pt}}
$5\mathrm{A}_{*o}$,
$3\mathrm{A}_{*o}\mathrm{A}_{n**o}$, 
$2\mathrm{A}_{*o}\mathrm{A}_{3**o}'$, 
$\mathrm{A}_{*o}\mathrm{A}_{n**o}\mathrm{A}_{m**o}$,
$\mathrm{A}_{n**o}\mathrm{A}_{3**o}'$,
$2\mathrm{A}_{*o}\mathrm{D}_{5*o}$,
$\mathrm{A}_{*o}\mathrm{A}_{5**o}'$,
$\mathrm{A}_{n**o}\mathrm{D}_{5*o}$.
\end{tabular}}
\end{multline}
\end{lemma}

\begin{proof}
Note that $\bar F$ is obtained from 
the minimal resolution $\mu\colon \tilde F\to F$ by 
contracting all the $-2$-curves:
\[
\mu\colon \tilde F\stackrel{\tilde{\mu}}{\longrightarrow}
\bar F\stackrel{\sigma_F}{\longrightarrow} F
\]

Moreover, $C_1$ and $C_2$ 
are exceptional curves of the induced map $\sigma_F\colon 
\bar F\to F$ 
and the coefficient of the proper 
transform of $C_1$ (resp. $C_2$) in the fundamental
cycle $Z$ is $2$ (resp. $1$). 
Write
\[
\frac12A_1+\frac14A_2=\tilde{\mu}^*
\left(\frac12C_1+\frac14C_2\right)+\sum_{j\neq 1,2}r_jA_j,
\]
where $A_1$ and $A_2$ are proper transforms of
$C_1$ and $C_2$, respectively, and
the $A_j$ for $j\neq 1,2$ are $-2$ curves on $\tilde F$.
The rational constants can be computed from the system of
linear equations:
\begin{equation}
\label{eq-r}
A_i\cdot\left(\frac12A_1+\frac14A_2\right)=
\sum_{j\neq 1,2}r_jA_i\cdot A_j,\qquad i\neq1,2.
\end{equation}
Since $K_{\tilde F}=\tilde{\mu}^*K_{\bar F}$, we have
$a\left(A_j,\frac12C_1+\frac14C_2\right)=-r_j$.
Hence, $(\bar F,\frac12C_1+\frac14C_2)$ is lc (resp. klt)
if and only if $\min_{j\neq 1,2}\{r_j\}\le 1$ (resp. $<1$).

\subsubsection*{}
\label{notation-graphs}
Now the assertion can be proved by elementary linear 
algebra computations. 
In our situation (when $d=3$ and $f_3$ has a multiple factor), 
the exceptional divisor 
is a tree of smooth rational curves. 
Therefore, all discrete invariants are
uniquely defined by weighted dual graphs.
Below the notation of \cite{L} are used.
Thus $\circ$ denotes the vertex corresponding to
exceptional curve $A_{\circ}$ with $A_{\circ}^2=-3$.
The vertices $\bullet$ always have $A_{\bullet}^2=-2$
while $A_{*}^2$ can be $-2$ or $-3$. We attach
$-3$ to corresponding $*$-vertices with 
$A_{*}^2=-3$ and omit $-2$ everywhere.
Note that each graph has at most one $\circ$-vertex.
By \cite{L} coefficients $z$ in the fundamental cycle $Z$
satisfy $z_*=1$ and $z_{\circ}=2$ or $3$. Therefore,
$C_2$ is a $*$-vertex and $C_1$ is $\circ$-vertex with 
$z_{\circ}=2$. 

Consider, for example, the case 
$2\mathrm{A}_{*o}\mathrm{A}_{3**o}'$. 
For \emph{weighted} dual graph of $\mu$ there
are two possibilities:
\begin{equation}
\label{eq--1}
\begin{array}{ccccccccc}
&&\stackrel{-3}{*}&\lin&\circ&\lin&*&&\\
&&&&|&&&&\\
*&\lin&\bullet&\lin&\bullet&\lin&\bullet&\lin&*\\
\end{array}
\end{equation}
or
\begin{equation}
\label{eq--2}
\begin{array}{ccccccccc}
&&{*}&\lin&\circ&\lin&*&&\\
&&&&|&&&&\\
*&\lin&\bullet&\lin&\bullet&\lin&\bullet&\lin&\stackunder{-3}{*}\\
\end{array}
\end{equation}
In case \eqref{eq--1}, contracting $-2$-curves, we 
get the surface $\bar F$ having two Du Val points
$P_1$ and $P_2$ of types $A_1$ and $A_5$, respectively.
Moreover, $P_1,P_2\in C_1$ and $P_1,P_2\notin C_2$.
It is clear that $(\bar F,C_1)$ is plt at $P_1$.
Hence $(\bar F,\frac12C_1+\frac14C_2)$ is klt at $P_1$.
The dual graph of the resolution $\tilde\mu$ over $P_2$
is as follows:
\begin{equation*}
\begin{array}{ccccccccc}
&&&&\stackrel{C_1}{\circ}&&&&\\
&&&&|&&&&\\
\stackunder{A_1}*&\lin&
\stackunder{A_2}\bullet&\lin&
\stackunder{A_3}\bullet&\lin&
\stackunder{A_4}\bullet&\lin&
\stackunder{A_5}*\\
\end{array}
\end{equation*}
Thus system \eqref{eq-r} has the form
\[
\left\{
\begin{array}{llllll}
-2r_1&+r_2&&&&=0\\
\phantom{+}r_1&-2r_2&+r_3&&&=0\\
&\phantom{+}r_2&-2r_3&+r_4&&=\frac12\\
&&\phantom{+}r_3&-2r_4&+r_5&=0\\
&&&\phantom{+}r_4&-2r_5&=0\\
\end{array}
\right.
\]
The solution is $(r_1,\dots,r_5)=
\left(-\frac14, -\frac12, -\frac34, -\frac12, -\frac14\right)$. Hence 
$(\bar F,\frac12C_1+\frac14C_2)$ is klt at $P_2$.

Similarly in case \eqref{eq--2} the surface $\bar F$ has 
three Du Val points
$P_1$, $P_2$, and $P_2$ of types $A_1$, $A_1$, and $A_4$, respectively.
Here $P_1,P_2\in C_1$, $P_1,P_2\notin C_2$ and $P_3=C_1\cap C_2$.
Clearly, $(\bar F,C_1+C_2)$ is plt outside of $P_3$.
Consider the dual graph of the resolution $\tilde\mu$ over $P_3$:
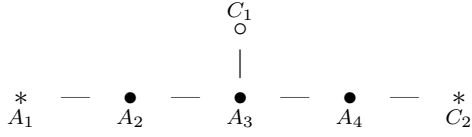
\begin{equation*}
\begin{array}{ccccccccc}
&&&&\stackrel{C_1}{\circ}&&&&\\
&&&&|&&&&\\
\stackunder{A_1}*&\lin&
\stackunder{A_2}\bullet&\lin&
\stackunder{A_3}\bullet&\lin&
\stackunder{A_4}\bullet&\lin&
\stackunder{C_2}*\\
\end{array}
\end{equation*}
This gives us the following form of system \eqref{eq-r}:
\[
\left\{
\begin{array}{llllll}
-2r_1&+r_2&&&&=0\\
\phantom{+}r_1&-2r_2&+r_3&&&=0\\
&\phantom{+}r_2&-2r_3&+r_4&&=\frac12\\
&&\phantom{+}r_3&-2r_4&&=\frac14\\
\end{array}
\right.
\]
The solution is $(r_1,\cdots,r_4)=
\left(-\frac14, -\frac12, -\frac34, -\frac12\right)$.
Hence 
$(\bar F,\frac12C_1+\frac14C_2)$ is klt at $P_3$.

Computations in all other cases are similar.
We omit them.
\end{proof}

\subsection*{}
By Lemma~\xref{l-exc} Theorem~\xref{main1}
is proved in all cases \eqref{mul-cacses}.
It is remain to consider cases 
$2\mathrm{A}_{*o}\mathrm{E}_{7o}$,
$\mathrm{A}_{*o}\mathrm{D}_{7*o}$,
$\mathrm{A}_{7**o}'$, 
$\mathrm{D}_{9*o}$, and
$\mathrm{A}_{n**o}\mathrm{E}_{7o}$. 
Here the proof can be done in the same style as above.
However, it is easier to use our computations in
Table~\xref{ta:res3} below.

\section{Proof of Theorem~\ref{main2}}
\label{proof-T2}
Let $\var\colon Y\to X$ be the blowup of $S$
(see \cite[17.10]{Ut}) and
let $F_Y$ be the proper transform of $F$.
By construction,
$Y$ is $\QQ$-factorial and the exceptional divisor of $\var$ is $S$.
Therefore, $\rho(Y/\CC^3)=1$ and $-S$ is $\var$-ample. We can write 
\[
f^*(K_{\CC^3}+cF)=K_Y+S+cF_Y.
\] 
Hence, $(Y,S+cF_Y)$
is plt (see \cite[3.10]{K}). If there is another blowup
$\var'\colon Y'\to X$ of $S$, then the composition map 
$Y\dasharrow Y'$ is an isomorphism in codimension one.
In this situation, $Y\dasharrow Y'$ must be an isomorphism.
This proves (i) and (ii). 
The assertion of (iii) follows by the Adjunction
\cite[17.6]{Ut}.

To prove (iv) we take $X'$ as  the normalization of $\CC^3$ in the 
finite extension $\KKK(\CC^3)[\sqrt[m]{f}]$,
where $f(x,y,z)=0$ is the equation of $F$ in $\CC^3$
(see e.g. \cite[\S 2]{Sh}).
Let $\pi\colon X'\to \CC^3$ be the natural projection.
Then $\pi$ is ramified only over $F$.
Hence $K_{X'}=\pi^*(K_{\CC^3}+cF)$ is linearly trivial.
By \cite[\S 2]{Sh}, $X'$ has lc non-klt singularities (of index $1$). 
Let $Y'$ be the normalization of a dominant component of 
$Y\times_{\CC^3}X'$. We obtain diagram \eqref{cd}, 
where the projection $\psi$ is a finite morphism while 
$\var'$ is birational. Put $S':=\psi^{-1}(S)_{\mt{red}}$.
Since
\begin{equation}
\label{eq-mul-d}
K_Y'+S'=\psi^*(K_Y+S+cF_Y),
\end{equation}
the pair $(Y',S')$ is plt (see \cite[\S 2]{Sh}).
On the other hand, the locus of log canonical
singularities of $(Y',S')$ is connected \cite[17.4]{Ut}
near any fiber of $\var'$. This shows that $S'$ is irreducible.
Therefore $S'$ is the only divisor with property 
$a(S',0)=-1$, i.e. $(Y,cF_Y)$ is exceptional.
Thus we proved that the pair $(Y',S')$ is plt of
index $1$, so it is canonical. 
By the Adjunction \cite[17.6]{Ut}, $S'$ has only Du Val 
singularities and $S'\sim 0$.

Finally, we have to show that the action of $\muu_m$ on $S'$
is faithful. 
Indeed, assume that some $0\neq \xi \in\muu_m$
acts trivially on $S'$. Let $F'_Y:=\psi^{-1}(F)_{\mt{red}}$ 
and let $P\in F'_Y\cap S'$ be the general point. 
Then $\xi(P)=P$. By construction, $\muu_m$ acts 
trivially on $F_Y'$. Thus it is sufficient to
show that the tangent spaces 
$\mathcal{T}_{S',P}$ and $\mathcal{T}_{F_Y',P}$
generate $\mathcal{T}_{Y',P}$. Assume the opposite.
Since $Y'$ has only terminal
singularities, $P\in Y'$ is smooth.
Hence $\mathcal{T}_{S',P}=\mathcal{T}_{F_Y',P}$
and $S'$ is tangent to $F_Y'$ along $S'\cap F_Y'$.
In particular, $(Y',S'+F_Y')$ is not lc along $S'\cap F_Y'$.
Similar to \eqref{eq-mul-d} we have
\[
K_Y'+S'+F_Y'=\psi^*(K_Y+S+F_Y).
\]
Therefore, $(Y,S+F_Y)$ also is not lc 
along $S'\cap F_Y'$ (see \cite[\S 2]{Sh}).
On the other hand, $(Y,S+cF_Y)$ is lc, where $c>1/2$.
Taking the general hyperplane section we derive a contradiction
(see e.g. \cite[Ex. 4.4.4]{Lect}).

\section{Tables}
\label{tabl}
Notation as in Theorem~\xref{main2}.
Below in Tables \xref{ta:res1}, \xref{ta:res2}, and \xref{ta:res3}, 
we enumerate our computations of log canonical thresholds,
blowup $\var$,
log Enriques surfaces and simple K3 singularities described in 
Theorem~\xref{main2} for all
sample equations \cite{L}. 
Below is the description of our tables.
\begin{itemize}
\item
The first column contains
types of minimal resolutions of singularities $(F\ni o)$.
In notation we follow \cite{L}. Note that symbols like
$\mathrm{A}_{1****}$ determine only the unweighted graph
of the minimal resolution $\mu\colon \tilde F\to F$.
For short, we omit the weights. 

\item
The second column contains two equations:
sample equations $f(x,y,z)$ of the 
singularity $F\subset \CC^3$ \cite{L} and $\ell(x,y,z)$
(see below).
\item
In the third column we indicate the log canonical threshold
$c=c(\CC^3,F)$.
\item
It turns out that in all cases the blowup $\var$ 
is a toric weighted blowup. The corresponding weights $w=(w_1,w_2,w_3)$
are listed in the forth column. In this case the exceptional
divisor $S$ must be a weighted projective plane $\PP(w_1,w_2,w_3)$
(see e.g. \cite{Flet}).
In the fifth column we write $S$ as $\PP(w_1',w_2',w_3')$,
where the triple $(w_1',w_2',w_3')$ is \emph{well-formed}
\cite{Flet}.
\item
Since $(Y,S)$ is a toric pair, $\Diff_S(0)$
is supported into the toric boundary of 
$S\simeq\PP(w_1',w_2',w_3')$. Moreover,
\[
\Diff_S(0)=\sum_{i=1}^3 \left(1-\frac1{m_i}\right)L_i,
\]
where $L_1$, $L_2$, $L_3$ are the coordinate ``lines'' 
$\{x=0\}$, $\{y=0\}$, $\{z=0\}$
and $m_i=\gcd(w_j,w_k)$, $\{i,j,k\}=\{1,2,3\}$.
Therefore, 
\[
\Diff_S(cF_Y)=cL+\sum_{i=1}^3 \delta_iL_i,
\] 
where $L$ is given by $\ell(x,y,z)=0$ 
and coefficients $\delta_i$ can be computed as
\[
\delta_i=
\begin{cases}
1-\frac1{m_i}+\frac{c}{m_i}&\text{if $F_Y\supset L_i$,}\\
1-\frac1{m_i}&\text{otherwise.}
\end{cases}
\]
(see \cite[Cor. 3.10]{Sh}).
\item
The last column contains notation in \cite{Y} 
for the simple K3 singularity $X'$
(whenever $c$ has the standard form).
\end{itemize}

\par\bigskip\bigskip\noindent
\setlongtables
\begin{longtable}
{p{55pt}p{100pt}p{10pt}p{45pt}p{50pt}
p{10pt}p{10pt}p{10pt}p{15pt}} 
\caption{Singularities with $Z^2=-1$}\label{ta:res1}\\
\\
\hline
\\
\multicolumn{1}{l}{\cite{L}}
&\multicolumn{1}{c}{$f$ and $\ell$}
&\multicolumn{1}{c}{$c$}
&\multicolumn{1}{c}{$w$}
&\multicolumn{1}{c}{$S$}
&\multicolumn{1}{c}{$\delta_1$}
&\multicolumn{1}{c}{$\delta_2$}
&\multicolumn{1}{c}{$\delta_3$}
&\cite{Y}\\
\\
\hline
\\
\endfirsthead
\multicolumn{8}{c}
{\small Table {\ref{ta:res1}} (continued)}\\ 
\hline
\\
\multicolumn{1}{l}{\cite{L}}
&\multicolumn{1}{c}{$f$ and $\ell$}
&\multicolumn{1}{c}{$c$}
&\multicolumn{1}{c}{$w$}
&\multicolumn{1}{c}{$S$}
&\multicolumn{1}{c}{$\delta_1$}
&\multicolumn{1}{c}{$\delta_2$}
&\multicolumn{1}{c}{$\delta_3$}
&\cite{Y}\\
\\
\hline
\\
\endhead
\hline 
\endlastfoot
\\
\endfoot

$\mathrm{Cu}$&$x^7+y^3+z^2$ & $\frac{41}{42}$ & $(6,14,21)$ & $\PP^2$ &
$\frac67$&$\frac23$&$\frac12$&\y{14}\\
\nopagebreak
&$x+y+z$ &&&&&&\vspace{6pt}\\

$\mathrm{Ta}$&$x^5y+y^3+z^2$& $\frac{29}{30}$& $(4,10,15)$& $\PP(2,1,3)$& 
$\frac45$&$-$&$\frac12$&\y{50}\\
\nopagebreak
&$xy+y^3+z$\vspace{6pt}\\

$\mathrm{Tr}$&$x^8+y^3+z^2$& $\frac{23}{24}$& $(3,8,12)$& $\PP(1,2,1)$& 
$\frac34$&$\frac23$&$-$&\y{13}\\ 
\nopagebreak
&$x^2+y+z^2$&&&&&&\vspace{6pt}\\

$\mathrm{A}_{1****}$&$x^9+y^3+z^2$&
$\frac{17}{18}$&$(2,6,9)$&$\PP(1,1,3)$& $\frac23$&$-$&$\frac12$&\y{12}\\
\nopagebreak
&$x^3+y^3+z$&&&&&&\vspace{6pt}\\

$\mathrm{A}_{n****}$
&\multicolumn{3}{l}{$(y+x^3)(y^2+x^{n+5})+z^2$}&&&\\ 
\nopagebreak
$\scriptstyle{n\ge 2}$&$(x+y)y^2+z$&$\frac{17}{18}$
&$(2,6,9)$
&$\PP(1,1,3)$
&$\frac23$&&$\frac12$&\y{12}\vspace{6pt}\\

$\mathrm{D}_{4***}$&$x^{10}+y^3+z^2$& $\frac{14}{15}$& $(3,10,15)$&
$\PP(1,2,1)$& $\frac45$&$\frac23$&$-$&\y{11}\\
\nopagebreak
&$x^2+y+z^2$&&&&&&\vspace{6pt}\\

$\mathrm{E}_{6**}$&$x^7y+y^3+z^2$& $\frac{13}{14}$& 
$(4,14,21)$&$\PP(2,1,3)$& 
$\frac67$&$-$&$\frac12$&\y{47}\\
\nopagebreak
&$xy+y^3+z$&&&&&&\vspace{6pt}\\

$\mathrm{E}_{8*}$&$x^{11}+y^3+z^2$ & $\frac{61}{66}$& $(6,22,33)$&$\PP^2$& 
$\frac{10}{11}$&$\frac23$&$\frac12$&\y{$-$}\\
\nopagebreak
&$x+y+z$&&&&&&\vspace{6pt}\\

\end{longtable}

\setlongtables
\begin{longtable}
{p{55pt}p{100pt}p{10pt}p{45pt}p{50pt}
p{10pt}p{10pt}p{10pt}p{15pt}} 
\caption{Singularities with $Z^2=-2$}\label{ta:res2}\\
\\
\hline
\\
\multicolumn{1}{l}{\cite{L}}
&\multicolumn{1}{c}{$f$ and $\ell$}
&\multicolumn{1}{c}{$c$}
&\multicolumn{1}{c}{$w$}
&\multicolumn{1}{c}{$S$}
&\multicolumn{1}{c}{$\delta_1$}
&\multicolumn{1}{c}{$\delta_2$}
&\multicolumn{1}{c}{$\delta_3$}
&\cite{Y}\\
\\
\hline
\\
\endfirsthead
\multicolumn{8}{c}
{\small Table {\ref{ta:res2}} (continued)}\\ 
\hline
\\
\multicolumn{1}{l}{\cite{L}}
&\multicolumn{1}{c}{$f$ and $\ell$}
&\multicolumn{1}{c}{$c$}
&\multicolumn{1}{c}{$w$}
&\multicolumn{1}{c}{$S$}
&\multicolumn{1}{c}{$\delta_1$}
&\multicolumn{1}{c}{$\delta_2$}
&\multicolumn{1}{c}{$\delta_3$}
&\cite{Y}\\
\\
\hline
\\
\endhead
\hline 
\endlastfoot
\\
\endfoot

\multicolumn{8}{c}{$\mathrm{I}$. Cases: $f_4(x,y)=xy^3$}\\
\\
\nopagebreak
$\mathrm{Cu}$&$z^2+x(x^4+y^3)$
& $\frac{29}{30}$
& $(6,8,15)$
& $\PP(1,4,5)$
&$-$& $\frac23$& $\frac12$&\y{38}\\
\nopagebreak
&$x^5+xy+z$&&&&&&\vspace{6pt}\\

$\mathrm{Ta}$&$z^2+x(x^3y+y^3)$
&$\frac{21}{22}$
&$(4,6,11)$&
$\PP(2,3,11)$
&$-$&$-$&$\frac12$&\y{78}\\
\nopagebreak
&$x^4y+xy^3+z$&&&&&&\vspace{6pt}\\

$\mathrm{Tr}$&$z^2+x(x^5+y^3)$ 
&$\frac{17}{18}$ 
&$(3,5,9)$ 
&$\PP(1,5,3)$ 
&$-$&$\frac23$&$-$&\y{39}\\ 
\nopagebreak
&$x^6+xy+z^2$&&&&&&\vspace{6pt}\\

$\mathrm{A}_{1****}$&$z^2+x(x^6+y^3)$
&$\frac{13}{14}$ 
&$(2,4,7)$ 
&$\PP(1,2,7)$ 
&$-$&$-$&$\frac12$&\y{40}\\ 
\nopagebreak
&$x^7+xy^3+z$&&&&&&\vspace{6pt}\\

$\mathrm{A}_{n****}$
&\multicolumn{4}{l}{$z^2+x(y+x^2)(y^2+x^{n+3})$}&&\\ 
\nopagebreak
$\scriptstyle{n\ge 2}$&$x(y+x^2)y^2+z$&$\frac{13}{14}$ 
&$(2,4,7)$ 
&$\PP(1,2,7)$ 
&$-$&$-$&$\frac12$&\y{40}\vspace{6pt}\\

$\mathrm{D}_{4***}$&$z^2+x(x^7+y^3)$ 
&$\frac{11}{12}$ 
&$(3,7,12)$ 
&$\PP(1,7,4)$ 
&$-$&$\frac23$&$-$&\y{41}\\ 
\nopagebreak
&$x^8+xy+z^2$&&&&&&\vspace{6pt}\\

$\mathrm{E}_{6**}$&$z^2+x(x^5y+y^3)$ 
&$\frac{31}{34}$ 
&$(4,10,17)$ 
&$\PP(2,5,17)$ 
&$-$&$-$&$\frac12$&\y{$-$}\\
\nopagebreak
&$x^6y+xy^3+z$&&&&&&\vspace{6pt}\\

$\mathrm{E}_{8*}$&$z^2+x(x^8+y^3)$ 
&$\frac{49}{54}$ 
&$(6,16,27)$ 
&$\PP(1,8,9)$ 
&$-$&$\frac23$& $\frac12$&\y{$-$}\\
\nopagebreak
&$x^9+xy+z$&&&&&&\vspace{6pt}\\
\hline

\\
\multicolumn{8}{c}{$\mathrm{II}$. Cases: $f_4(x,y)=y^4$}\\
\nopagebreak
\\
\nopagebreak
$\mathrm{Ta}$&$z^2+x^5+y^4$ 
&$\frac{19}{20}$ 
&$(4,5,10)$ 
&$\PP(2,1,1)$ 
&$\frac45$&$\frac12$& $-$&\y{9}\\
\nopagebreak
&$x+y^2+z^2$&&&&&&\vspace{6pt}\\

$\mathrm{Tr}$&$z^2+x^4y+y^4$ 
&$\frac{15}{16}$ 
&$(3,4,8)$ 
&$\PP(3,1,2)$ 
&$\frac34$&$-$ &$-$&\y{37}\\
\nopagebreak
&$xy+y^4+z^2$&\vspace{6pt}\\

$\mathrm{A}_{1****}$&$z^2+x^6+y^4$ 
&$\frac{11}{12}$ 
&$(2,3,6)$ 
&$\PP^2$ 
&$\frac23$&$\frac12$& $-$&\y{8}\\
\nopagebreak
&$x^2+y^2+z^2$&&&&&&\vspace{6pt}\\

$\mathrm{A}_{n****}$, $\scriptstyle{n\ge 2}$
&\multicolumn{3}{l}{$z^2+(y^2+x^3)(y^2+x^{n+2})$}&
$\PP^2$
&\\ \nopagebreak
&${y}^{2}+xy+{z}^{2}$& 
$\frac{11}{12}$&$(2,3,6)$& 
&$\frac23$&$\frac12$& $-$&\y{8}\vspace{6pt}\\

$\mathrm{A}_{n****}$, $\scriptstyle{n\ge 2}$
&\multicolumn{6}{l}{$z^2+(y^2+x^3)^2+x^ay^b$,\ 
${2a+3b=n+11}$}\\ 
\nopagebreak
&${y}^{2}+xy+{x}^{2}+{z}^{2}$&$\frac{11}{12}$ 
&$(2,3,6)$ 
&$\PP^2$ 
&$\frac23$&$\frac12$& $-$&\y{8}\vspace{6pt}\\

$\mathrm{D}_{4**}$&$z^2+x^5y+y^4$
&$\frac9{10}$ 
&$(3,5,10)$ 
&$\PP(3,1,2)$ 
&$\frac45$&$-$&$-$&\y{36}\\
\nopagebreak
&$xy+y^4+z^2$&&&&&&\vspace{6pt}\\

$\mathrm{E}_{6**}$&$z^2+x^7+y^4$ 
&$\frac{25}{28}$
&$(4,7,14)$ 
&$\PP(2,1,1)$ 
&$\frac67$&$\frac12$& $-$&\y{$-$}\\
\nopagebreak
&$x+y^2+z^2$&&&&&&\vspace{6pt}\\

\hline
\\
\multicolumn{8}{c}{$\mathrm{III}$. Cases: $f_4(x,y)=0$}\\
\nopagebreak
\\
\nopagebreak
$5\mathrm{A}_{*o}$&$z^2+x^5+y^5$ 
&$\frac9{10}$ 
&$(2,2,5)$ 
&$\PP(1,1,5)$ 
&$-$&$-$&$\frac12$&\y{6}\\
\nopagebreak
&$x^5+y^5+z$&&&&&&\vspace{6pt}\\

$3\mathrm{A}_{*o}\mathrm{A}_{n**o}$
&\multicolumn{2}{l}{$z^2+(y^3+x^3)(x^2+y^{n+2})$} 
&$(2,2,5)$ 
&$\PP(1,1,5)$ 
&$-$&$-$&$\frac12$&\y{6}\\
\nopagebreak 
&${y}^{3}{x}^{2}+{x}^{5}+z$&$\frac{9}{10}$&\vspace{6pt}\\

$2\mathrm{A}_{*o}\mathrm{A}_{3**o}'$
&$z^2+(y^2+x^2)(x^3+y^4)$ 
&$\frac8{9}$ 
&$(4, 3, 9)$ 
&$\PP(4, 1, 3)$ 
&$\frac23$&$-$&$-$&\y{33}\\
\nopagebreak
&$x{y}^{2}+{y}^{6}+{z}^{2}$&&&&&&\vspace{6pt}\\

$2\mathrm{A}_{*o}\mathrm{D}_{5*o}$
&\multicolumn{2}{l}{$z^2+(y^2+x^2)(x^3+xy^3)$} 
&$(6, 4, 13)$ 
&$\PP(3, 2, 13)$ 
&$-$&$-$&$\frac12$&\y{$-$}\\
\nopagebreak
&${y}^{2}{x}^{3}+{y}^{5}x+z$&$\frac{23}{26}$ 
&&&&&\vspace{6pt}\\

$2\mathrm{A}_{*o}\mathrm{E}_{7o}$
&$z^2+(y^2+x^2)(x^3+y^5)$ 
&$\frac{37}{42}$ 
&$(10, 6, 21)$ 
&$\PP(5, 1, 7)$ 
&$\frac23$&$-$&$\frac12$&\y{$-$}\\
\nopagebreak
&$x{y}^{2}+{y}^{7}+z$&&&&&&\vspace{6pt}\\

\multicolumn{2}{l}
{$\mathrm{A}_{*o}\mathrm{A}_{n**o}\mathrm{A}_{m**o}$}
&\\ 
\nopagebreak
&\multicolumn{3}{l}{$z^2+(y+x)(y^2+x^{n+2})(x^2+y^{m+2})$}\\ 
\nopagebreak
&${y}^{3}{x}^{2}+{x}^{3}{y}^{2}+z$&$\frac{9}{10}$ 
&$(2, 2, 5)$ 
&$\PP(1, 1, 5)$ 
&$-$&$-$&$\frac12$&\y{6}\vspace{6pt}\\

$\mathrm{A}_{*o}\mathrm{A}_{5**o}'$&$z^2+x^4y+y^6$ 
&$\frac78$ 
&$(5,4,12)$ 
&$\PP(5,1,3)$ 
&$\frac34$&$-$&$-$&\y{31}\\ 
\nopagebreak
&$xy+y^6+z^2$&&&&&&\vspace{6pt}\\

$\mathrm{A}_{*o}\mathrm{D}_{7*o}$&$z^2+x^4y+xy^5$ 
&$\frac{33}{38}$ 
&$(8,6,19)$ 
&$\PP(4,3,19)$ 
&$-$&$-$&$\frac12$&\y{$-$}\\
\nopagebreak
&$x^4y+xy^5+z$&&&&&&\vspace{6pt}\\

$\mathrm{A}_{7**o}'$&$z^2+x^5+y^6$ 
&$\frac{13}{15}$ 
&$(6,5,15)$ 
&$\PP(2,1,1)$ 
&$\frac45$&$\frac23$& $-$&\y{$-$} \\
\nopagebreak
&$x+y^2+z^2$&&&&&&\vspace{6pt}\\

$\mathrm{D}_{9*o}$&$z^2+x^5+xy^5$ 
&$\frac{43}{50}$ 
&$(10,8,25)$ 
&$\PP(1,4,5)$ 
&$-$&$\frac45$& $\frac12$&\y{$-$}\\
\nopagebreak
&$x^5+xy+z$&&&&&&\vspace{6pt}\\

$\mathrm{A}_{n**o}\mathrm{A}_{3**o}'$
&\multicolumn{6}{l}{$z^2+(y^2+x^{n+2})(x^3+y^4)$}\\ 
\nopagebreak
&$x{y}^{2}+{y}^{6}+{z}^{2}$&$\frac{8}{9}$ 
&$(4,3,9)$ 
&$\PP(4, 1, 3)$ 
&$\frac23$&$-$&$-$&\y{33}\vspace{6pt}\\

$\mathrm{A}_{n**o}\mathrm{D}_{5*o}$
&\multicolumn{6}{l}{$z^2+(y^2+x^{n+2})(x^3+xy^3)$}\\ 
\nopagebreak
&${y}^{2}{x}^{3}+{y}^{5}x+z$&$\frac{23}{26}$ 
&$(6, 4, 13)$ 
&$\PP(3, 2, 13)$ 
&$-$&$-$&$\frac12$&\y{$-$}\vspace{6pt}\\

$\mathrm{A}_{n**o}\mathrm{E}_{7o}'$
&\multicolumn{6}{l}{$z^2+(y^2+x^{n+2})(x^3+y^5)$}\\ 
\nopagebreak
&$x{y}^{2}+{y}^{7}+z$&$\frac{37}{42}$ 
&$(10, 6, 21)$ 
&$\PP(5, 1, 7)$ 
&$\frac23$&$-$&$\frac12$&\y{$-$}\vspace{6pt}\\
\end{longtable}

\par\bigskip\bigskip\noindent
\setlongtables
\begin{longtable}
{p{55pt}p{100pt}p{10pt}p{45pt}p{50pt}
p{10pt}p{10pt}p{10pt}p{15pt}} 
\caption{Singularities with $Z^2=-3$}\label{ta:res3}\\
\\
\hline
\\
\multicolumn{1}{l}{\cite{L}}
&\multicolumn{1}{c}{$f$ and $\ell$}
&\multicolumn{1}{c}{$c$}
&\multicolumn{1}{c}{$w$}
&\multicolumn{1}{c}{$S$}
&\multicolumn{1}{c}{$\delta_1$}
&\multicolumn{1}{c}{$\delta_2$}
&\multicolumn{1}{c}{$\delta_3$}
&\cite{Y}\\
\\
\hline
\\
\endfirsthead
\multicolumn{8}{c}
{\small Table {\ref{ta:res3}} (continued)}\\ 
\hline
\\
\multicolumn{1}{l}{\cite{L}}
&\multicolumn{1}{c}{$f$ and $\ell$}
&\multicolumn{1}{c}{$c$}
&\multicolumn{1}{c}{$w$}
&\multicolumn{1}{c}{$S$}
&\multicolumn{1}{c}{$\delta_1$}
&\multicolumn{1}{c}{$\delta_2$}
&\multicolumn{1}{c}{$\delta_3$}
&\cite{Y}\\
\\
\hline
\\
\endhead
\hline 
\endlastfoot
\\
\endfoot

\multicolumn{8}{c}{$\mathrm{I}$. Cases: $f_3(x,y,z)=y^3+xz^2$}\\
\nopagebreak
\\
\nopagebreak
$\mathrm{Cu}$&$x^4+y^3+xz^2$ 
&$\frac{23}{24}$ 
&$(6,8,9)$ 
&$\PP(1,4,3)$ 
&$-$&$\frac23$& $\frac12$&\y{20} \\
\nopagebreak
&$x^4+y+xz$&&&&&&\vspace{6pt}\\

$\mathrm{Ta}$&$x^3y+y^3+xz^2$ 
&$\frac{17}{18}$ 
&$(4,6,7)$ 
&$\PP(2,3,7)$
&$-$&$-$&$\frac12$&\y{60} \\
\nopagebreak
&$x^3y+y^3+xz$&&&&&&\vspace{6pt}\\

$\mathrm{Tr}$&$x^5+y^3+xz^2$ 
&$\frac{14}{15}$ 
&$(3,5,6)$ 
&$\PP(1,5,2)$ 
&$-$&$\frac23$&$-$&\y{22}\\ 
\nopagebreak
&$x^5+y+xz^2$&&&&&&\vspace{6pt}\\

$\mathrm{A}_{1****}$&$x^6+y^3+xz^2$ 
&$\frac{11}{12}$ 
&$(2,4,5)$ 
&$\PP(1,2,5)$ 
&$-$&$-$&$\frac12$&\y{24} \\
\nopagebreak
&$x^6+y^3+xz$&&&&&&\vspace{6pt}\\

$\mathrm{A}_{n****}$, $\scriptstyle{n\ge 2}$
&\multicolumn{3}{l}
{$xz^2+y^3+x^2y^2+(y+x^2)x^{n+3}$}&&&&&\\
\nopagebreak
&${x}^{2}{y}^{2}+{y}^{3}+xz
$&$\frac{11}{12}$ 
&$(2,4,5)$ 
&$\PP(1,2,5)$ 
&$-$&$-$&$\frac12$&\y{24}\vspace{6pt}\\

$\mathrm{D}_{4***}$&$x^7+y^3+xz^2$ 
&$\frac{19}{21}$ 
&$(3,7,9)$ 
&$\PP(1,7,3)$ 
&$-$&$\frac23$&$-$&\y{$-$}\\
\nopagebreak
&$x^7+y+xz^2$&&&&&&\\

$\mathrm{E}_{6**}$&$x^5y+y^3+xz^2$ 
&$\frac9{10}$ 
&$(4,10,13)$ 
&$\PP(2,5,13)$ 
&$-$&$-$&$\frac12$&\y{68} \\
\nopagebreak
&$x^5y+y^3+xz$&&&&&&\vspace{6pt}\\

$\mathrm{E}_{8*}$&$x^8+y^3+xz^2$ 
&$\frac{43}{48}$ 
&$(6,16,21)$ 
&$\PP(1,8,7)$ 
&$-$&$\frac23$&$\frac12$&\y{$-$}\\
\nopagebreak
&$x^8+y+xz$&&&&&&\vspace{6pt}\\
\hline

\\
\multicolumn{8}{c}
{$\mathrm{II}$. Cases: $f_3(x,y,z)=(y^2+xz)z$}\\
\nopagebreak
\\
\nopagebreak
$\mathrm{Ta}$&$x^4+y^2z+xz^2$ 
&$\frac{15}{16}$ 
&$(4,5,6)$ 
&$\PP(2,5,3)$ 
&$-$&$\frac12$&$-$&\y{58} \\
\nopagebreak
&$x^4+yz+xz^2$&&&&&&\vspace{6pt}\\

$\mathrm{Tr}$&$x^3y+y^2z+xz^2$ 
&$\frac{12}{13}$ 
&$(3,4,5)$ 
&$\PP(3,4,5)$ 
&$-$&$-$&$-$&\y{87} \\
\nopagebreak
&$x^3y+y^2z+xz^2$&&&&&&\vspace{6pt}\\

$\mathrm{A}_{1****}$&$x^5+y^2z+xz^2$ 
&$\frac9{10}$ 
&$(2,3,4)$ 
&$\PP(1,3,2)$ 
&$-$&$\frac12$&$-$&\y{63} \\
\nopagebreak
&$x^5+yz+xz^2$&&&&&&\vspace{6pt}\\

$\mathrm{A}_{n****}$, $\scriptstyle{n\ge 2}$
&\multicolumn{5}{l}
{$z(xz-2y^2)+x^2y^2+(y^2+x^3)x^{n+1}$}
&$ $&$ $&\\
\nopagebreak
&$x{z}^{2}+yz+{x}^{2}y$&$\frac9{10}$ 
&$(2,3,4)$ 
&$\PP(1,3,2)$ 
&$-$&$\frac12$&$-$&\y{63}\vspace{6pt}\\

$\mathrm{A}_{n****}$, $\scriptstyle{n\ge 2}$
&\multicolumn{6}{l}
{$z(xz-2y^2)+2x^2y^2+x^5+x^ay^b$,\ $2a+3b=n+9$}
&$ $&\\
\nopagebreak
&${x}^{5}+x{z}^{2}+yz+{x}^{2}y$&$\frac9{10}$ 
&$(2,3,4)$ 
&$\PP(1,3,2)$ 
&$-$&$\frac12$&$-$&\y{63}\vspace{6pt}\\

$\mathrm{D}_{4***}$&$x^4y+y^2z+xz^2$ 
&$\frac{15}{17}$ 
&$(3,5,7)$ 
&$\PP(3,5,7)$ 
&$-$&$-$&$-$&\y{$-$}\\ 
\nopagebreak
&$x^4y+y^2z+xz^2$&&&&&&\vspace{6pt}\\

$\mathrm{E}_{6**}$&$x^6+y^2z+xz^2$ 
&$\frac78$ 
&$(4,7,10)$ 
&$\PP(2,7,5)$ 
&$-$&$\frac12$&$-$&\y{64}\\
\nopagebreak
&$x^6+yz+xz^2$&&&&&&\vspace{6pt}\\
\hline
\\
\multicolumn{8}{c}{$\mathrm{III}$. Cases: $f_3(x,y,z)=y^3+z^3$}\\
\nopagebreak
\\
\nopagebreak
$\mathrm{Tr}$&$x^4+y^3+z^3$ 
&$\frac{11}{12}$ 
&$(3,4,4)$ 
&$\PP(3,1,1)$ 
&$\frac34$&$-$&$-$&\y{4} \\
\nopagebreak
&$x+y^3+z^3$&&&&&&\vspace{6pt}\\

$\mathrm{A}_{1****}$&$x^3y+y^3+z^3$ 
&$\frac89$ 
&$(2,3,3)$ 
&$\PP(2,1,1)$ 
&$\frac23$&$-$&$-$&\y{18} \\
\nopagebreak
&$xy+y^3+z^3$&&&&&&\vspace{6pt}\\

$\mathrm{A}_{n****}$, $\scriptstyle{n\ge 2}$
&\multicolumn{6}{l}
{$z^3+y^3+x^3(y+z)+x^ay^b$,\ $2a+3b=n+8$}
&&\\
\nopagebreak
&\multicolumn{3}{l}{$\left (z+y\right )\left 
({y}^{2}-yz+x+{z}^{2}\right )
$}\\
&&$\frac89$ 
&$(2,3,3)$ 
&$\PP(2,1,1)$ 
&$\frac23$&$-$&$-$&\y{18}
\vspace{6pt}\\

$\mathrm{D}_{4***}$&$x^5+y^3+z^3$ 
&$\frac{13}{15}$ 
&$(3,5,5)$ 
&$\PP(3,1,1)$ 
&$\frac45$&$-$&$-$&\y{$-$} \\
\nopagebreak
&$x+y^3+z^3$&&&&&&\vspace{6pt}\\

\hline
\\
\multicolumn{8}{c}{$\mathrm{IV}$. Cases: $f_3(x,y,z)=yz^2$}\\
\nopagebreak
\\
\nopagebreak
$5\mathrm{A}_{*o}$&$x^4+y^4+yz^2$ 
&$\frac78$ 
&$(2,2,3)$ 
&$\PP(1,1,3)$ 
&$-$&$-$&$\frac12$&\y{19} \\
\nopagebreak
&$x^4+y^4+yz$&&&&&&\vspace{6pt}\\

$3\mathrm{A}_{*o}\mathrm{A}_{n**o}$
&$yz^2+x^4+x^2y^2+y^{n+4}$ 
&$\frac{7}{8}$ 
&$(2, 2, 3)$ 
&$\PP(1,1,3)$ 
&$-$&$-$&$\frac12$&\y{19}\\
\nopagebreak
&$yz+{x}^{2}{y}^{2}+{x}^{4}$&&&&&&\vspace{6pt}\\

$3\mathrm{A}_{*o}\mathrm{A}_{n**o}$
&\multicolumn{3}{l}{$yz^2+y^4+x^3y+x^az^b$,\ $2a+3b=n+8$}
&$\PP(1,1,3)$ 
&$-$&$\frac78$&$\frac12$&\y{19} \\
\nopagebreak
&$z+{y}^3+{x}^{3}$&$\frac78$&$(2,2,3)$&&&&\vspace{6pt}\\

$2\mathrm{A}_{*o}\mathrm{A}_{3**o}'$&$yz^2+x^4+x^3y+y^5$ 
&$\frac{13}{15}$ 
&$(4, 3, 6)$ 
&$\PP(2,1,1)$ 
&$\frac23$&$\frac{14}{15}$&$-$&\y{$-$} \\
\nopagebreak
&${z}^{2}+{y}^2+x$&&&&&&\vspace{6pt}\\

$2\mathrm{A}_{*o}\mathrm{A}_{3**o}'$
&$yz^2+y^4+x^2y^2+x^3z$ 
&$\frac{6}{7}$ 
&$(3, 4, 5)$ 
&$\PP(3,4,5)$ 
&$-$&$-$&$-$&\y{85} \\
\nopagebreak
&${x}^{2}{y}^{2}+{x}^{3}z+y{z}^{2}$&&&&&&\vspace{6pt}\\

$2\mathrm{A}_{*o}\mathrm{D}_{5*o}$
&$yz^2+x^4+x^3y+xy^4$ 
&$\frac{19}{22}$ 
&$(6, 4, 9)$ 
&$\PP(1,2,3)$ 
&$-$&$\frac{21}{22}$&$\frac12$&\y{$-$} \\
\nopagebreak
&$x^3+xy+z$&&&&&&\vspace{6pt}\\

$2\mathrm{A}_{*o}\mathrm{D}_{5*o}$
&$yz^2+y^4+x^2y^2+x^5$ 
&$\frac{17}{20}$ 
&$(4, 6, 7)$ 
&$\PP(2,3,7)$ 
&$-$&$-$&$\frac12$&\y{$-$} \\
\nopagebreak
&$yz+{x}^{2}{y}^{2}+{x}^{5}$&&&&&&\vspace{6pt}\\

$2\mathrm{A}_{*o}\mathrm{E}_{7*}$
&$yz^2+x^4+x^3y+y^6$ 
&$\frac{31}{36}$ 
&$(10, 6, 15)$ 
&$\PP^2$ 
&$\frac23$&$\frac{35}{36}$&$\frac12$&\y{$-$} \\
\nopagebreak
&$x+y+z$&&&&&&\vspace{6pt}\\

\multicolumn{2}{l}
{$\mathrm{A}_{*o}\mathrm{A}_{n**o}\mathrm{A}_{m**o}$}\\
&\multicolumn{3}{l}{$(x+y)z^2+x^2y^2+x^{n+4}+y^{m+4}$}
& 
&$-$&$-$&$\frac12$&\y{19} \\
\nopagebreak
&${x}^{2}{y}^{2}+xz+yz$&$\frac78$&$(2,2,3)$
&$\PP(1,1,3)$&&&\vspace{6pt}\\

\multicolumn{2}{l}
{$\mathrm{A}_{*o}\mathrm{A}_{n**o}\mathrm{A}_{m**o}$}\\
\nopagebreak
&\multicolumn{4}{l}{$yz^2+x^3y+x^2y^2+y^{n+4}+x^az^b$,\ $2a+3b=m+8$} 
&&&& \\
\nopagebreak
&${x}^{3}y+{x}^{2}{y}^{2}+yz$
&$\frac78$&$(2,2,3)$&$\PP(1,1,3)$&$-$&$-$&$\frac12$
&\y{19}\vspace{6pt}\\

$\mathrm{A}_{*o}\mathrm{A}_{5**o}'$&$x^4+y^5+yz^2$ 
&$\frac{17}{20}$ 
&$(5,4,8)$ 
&$\PP(5,1,2)$ 
&$\frac34$&$-$&$-$&\y{$-$} \\
\nopagebreak
&$x+y^5+yz^2$&&&&&&\vspace{6pt}\\

$\mathrm{A}_{*o}\mathrm{A}_{5**o}'$&$x^3z+xy^3+yz^2$ 
&$\frac{16}{19}$ 
&$(4,5,7)$ 
&$\PP(4,5,7)$ 
&$-$&$-$&$-$&\y{$-$} \\
\nopagebreak
&$x^3z+xy^3+yz^2$&&&&&&\vspace{6pt}\\

$\mathrm{A}_{*o}\mathrm{D}_{7*o}$&$x^4+xy^4+yz^2$ 
&$\frac{27}{32}$
&$(8,6,13)$ 
&$\PP(4,3,13)$ 
&$-$&$-$&$\frac12$ &\y{$-$}\\
\nopagebreak
&$x^4+xy^4+yz$&&&&&&\vspace{6pt}\\

$\mathrm{A}_{*o}\mathrm{D}_{7*o}$&$x^5+xy^3+yz^2$ 
&$\frac56$ 
&$(6,8,11)$ 
&$\PP(3,4,11)$ 
&$-$&$-$&$\frac12$&\y{56} \\
\nopagebreak
&$x^5+xy^3+yz$&&&&&&\vspace{6pt}\\

$\mathrm{A}_{7**o}'$&$x^3z+y^4+yz^2$ 
&$\frac56$ 
&$(5,6,9)$ 
&$\PP(5,2,3)$ 
&$\frac23$&$-$&$-$&\y{57} \\
\nopagebreak
&$xz+y^4+yz^2$&&&&&&\vspace{6pt}\\

$\mathrm{D}_{9*o}$&$x^5+y^4+yz^2$ 
&$\frac{33}{40}$ 
&$(8,10,15)$ 
&$\PP(4,1,3)$ 
&$\frac45$&$-$&$\frac12$&\y{$-$} \\
\nopagebreak
&$x+y^4+yz$&&&&&&\vspace{6pt}\\

$\mathrm{A}_{n**o}\mathrm{A}_{3**o}'$
&\multicolumn{6}{l}
{$yz^2+x^3y+x^az^b+y^5$,\ $2a+3b=n+8$}
&$ $&\\
\nopagebreak
&$x+y^2+z^2$&$\frac{13}{15}$&$(4,3,6)$&$\PP(2,1,1)$
&$\frac23$&$\frac{14}{15}$&$-$&\y{$-$}\vspace{6pt}\\

$\mathrm{A}_{n**o}\mathrm{A}_{3**o}'$
&\multicolumn{6}{l}
{$yz^2+x^2y^2+y^{n+4}+x^3z$}
&$ $&\\
\nopagebreak
&${x}^{2}{y}^{2}+{x}^{3}z+y{z}^{2}$&
$\frac67$&$(3,4,5)$&$\PP(3,4,5)$&$-$&$-$&$-$&\y{85}
\vspace{6pt}\\

$\mathrm{A}_{n**o}\mathrm{D}_{5*o}$
&\multicolumn{6}{l}
{$yz^2+x^3y+x^az^b+xy^4$,\ $2a+3b=n+8$}
&$ $&\\
\nopagebreak
&$x^3+xy+z$&$\frac{19}{22}$&$(6,4,9)$
&$\PP(1, 2, 3)$&$-$&$\frac{21}{22}$
&$\frac12$&\y{$-$}\vspace{6pt}\\

$\mathrm{A}_{n**o}\mathrm{D}_{5*o}$
&\multicolumn{6}{l}{$yz^2+x^2y^2+y^{n+4}+x^5$}&&\\
\nopagebreak
&$x^5+x^2y^2+yz$&$\frac{17}{20}$&$(4,6,7)$&$\PP(2,3,7)$&$-$&$-$&$\frac12$
&\y{$-$}\vspace{6pt}\\

$\mathrm{A}_{n**o}\mathrm{E}_{7o}$
&\multicolumn{6}{l}
{$yz^2+x^3y+x^az^b+y^6$,\ $2a+3b=n+8$}
&$ $&\\
\nopagebreak
&$x+y+z$&$\frac{31}{36}$&$(10,6,15)$&$\PP^2$&
$\frac23$&$\frac{35}{36}$&$\frac12$&\y{$-$}\vspace{6pt}\\

\hline
\\
\multicolumn{8}{c}{$\mathrm{V}$. Cases: $f_3(x,y,z)=z^3$}\\
\nopagebreak
\\
\nopagebreak
$4\mathrm{A}_{1*o}$&$x^4+y^4+z^3$ 
&$\frac56$ 
&$(3,3,4)$ 
&$\PP(1,1,4)$
&$-$&$-$&$\frac23$&\y{2} \\
\nopagebreak
&$x^4+y^4+z$&&&&&&\vspace{6pt}\\

$\mathrm{A}_{1*o}\mathrm{A}_{7*o}$&$x^3z+xy^3+z^3$
&$\frac{22}{27}$ 
&$(6,7,9)$ 
&$\PP(2,7,3)$ 
&$-$&$\frac23$&$-$&\y{$-$} \\
\nopagebreak
&$x^3z+xy+z^3$&&&&&&\vspace{6pt}\\

$\mathrm{A}_{10*o}$&$x^3z+y^4+z^3$ 
&$\frac{29}{36}$ 
&$(8,9,12)$ 
&$\PP(2,3,1)$ 
&$\frac23$&$\frac34$&$-$ &\y{$-$} \\
\nopagebreak
&$xz+y+z^3$&&&&&&\vspace{6pt}\\
\end{longtable}

\subsection*{Comments on computations}
Note that in all cases the weight $w$ and the constant $c$ 
are taken so that
\[
c(F)= \frac{\sum w_i}{\ord_{w,0}(f)},
\]
where $\ord_{w,0}(f)$ is the weighted multiplicity of $f(x,y,z)$ 
at $0$. Therefore, $a(S,cF)=-1$. 
Case by case one can show that the pair $(S,\Diff_S(cF)$ is klt.
Then, automatically, $\var\colon Y\to \CC^3$ in Theorem~\xref{main2}
coincides with our $w$-weighted blowup.

An alternative method is as follows.
In most cases the weighted homogeneous weighted term $f_w(x,y,z)$
of $f(x,y,z)$ with respect to $w$ defines a singularity which is 
isolated (or log canonical outside the origin). In these cases
we can use the following fact:

\begin{proposition}[{\cite[Prop. 8.14]{K}}]
\label{prop-K}
Let $f(x_1,\dots,x_n)$ be a holomorphic function near
$0\in\CC^n$ and $F:=\{f(x_1,\dots,x_n)=0\}$. Let $w=(w_1,\dots,w_n)$
be a weight. Then
\begin{equation}
\label{eq-eq-eq}
c(F)\le \frac{\sum w_i}{\ord_{w,0}(f)},
\end{equation}
Moreover, if the singularity $\{f_w(x_1,\dots,x_n)=0\}$
is lc outside the origin, then equality holds.
\end{proposition}

If in \eqref{eq-eq-eq} equality holds, then 
$\var\colon Y\to \CC^3$ in Theorem~\xref{main2} 
can be taken as the $w$-weighted blowup.

\end{document}